\documentclass{article}
\usepackage{amssymb,amsmath,amsthm}
\usepackage{setspace}

\newtheorem{thm}{Theorem}[section]
\newtheorem{prp}{Proposition}[section]
\newtheorem{cor}{Corollary}[section]
\newtheorem{lem}{Lemma}[section]
\theoremstyle{definition}
\newtheorem{dfn}{Definition}[section]
\numberwithin{equation}{section}
\newcounter{eqn}

\def\up#1{\uppercase{#1}}
\def\E{\expandafter\up}
\def\a{averag}

\def\ct{constant}
\def\df{difference}
\def\eq{equation}
\def\f{function}
\def\g{gradient}
\def\HL{Hardy--Littlewood}
\def\iq{inequalit}
\def\mms{metric measure space}
\def\P{Poincar\'e}
\def\pw{pointwise}
\def\q{quotient}
\def\rh{right-hand side}
\def\Sm{Sobolev metric }
\def\Ss{Sobolev space}
\def\wrt{with respect to }
\let\iy\infty
\let\n\nabla
\let\O\Omega
\let\S\Sigma
\let\z\zeta
\def\fint#1#2{\mathchoice
{\hbox{\,\ooalign{\hfil$\relbar$\hfil\crcr$\displaystyle\intop_{#1}$}\,}}
{\mathop{\hbox{\ooalign{\hfil$\relbar$\hfil\crcr$\textstyle\intop$}\!}}\nolimits_{#1}}
{\mathop{\hbox{\ooalign{\hfil$\relbar$\hfil\crcr$\textstyle\intop$}\!}}\nolimits_{#1}}
{\mathop{\hbox{\ooalign{\hfil$\relbar$\hfil\crcr$\textstyle\intop$}\!}}\nolimits_{#1}}}

\def\kint#1#2{\mathchoice
{\hbox{\,\ooalign{\hfil$\relbar$\hfil\crcr$\displaystyle\intop_{#1}^{#2}$}\,}}
{\mathop{\hbox{\ooalign{\hfil$\relbar$\hfil\crcr$\textstyle\intop$}\!}}\nolimits_{#1}^{#2}}
{\mathop{\hbox{\ooalign{\hfil$\relbar$\hfil\crcr$\textstyle\intop$}\!}}\nolimits_{#1}^{#2}}
{\mathop{\hbox{\ooalign{\hfil$\relbar$\hfil\crcr$\textstyle\intop$}\!}}\nolimits_{#1}^{#2}}}

\def\diam{\mathop{\rm diam}}
\def\M{\mathcal M}

\def\R{\mathbb R}
\def\C{\mathbb C}
\def\T{\mathfrak T}
\def\sn#1{|\kern-0.5pt\|#1\|\kern-0.5pt|}
\def\loc{_{\rm loc}}

\let\ti\textit
\def\er#1{\eqref{e#1}}

\def\bal#1\e{\begin{align}#1\end{align}}
\def\bga#1\e{\begin{gather}#1\end{gather}}
\def\beq#1\e{\begin{equation}#1\end{equation}}
\def\bml#1\e{\begin{multline}#1\end{multline}}
\def\bth{\begin{thm}}\def\eth{\end{thm}}
\def\bpr{\begin{prp}}\def\epr{\end{prp}}
\def\bco{\begin{cor}}\def\eco{\end{cor}}
\def\blm{\begin{lem}}\def\elm{\end{lem}}
\def\bdf{\begin{dfn}}\def\edf{\end{dfn}}
\def\ben{\begin{enumerate}}\def\een{\end{enumerate}}

\title{Sobolev spaces and averaging I\footnote{The main concepts of this paper were presented in
the talk at H.~Triebel's seminar {\it Function spaces}, Jena, April~12, 2013.}}
\author{B. Bojarski\\
Institute of Mathematics PAS, Warsaw}

\begin{document}

\maketitle

\rightline{\it Dedicated to Galina Fedorovna Bojarska}

\renewcommand{\thefootnote}{\arabic{footnote})}
\setcounter{footnote}{0}

\vskip15pt
\begin{abstract}
An apparently new concept of maximal mean \df\ \q\ is defined for \f s in the
Lebesgue space $L\loc(\R^n)$. Our definitions are meaningful for vector
valued \f s \cite{HeKST} on general measure metric spaces
as well and seem to lead to the most natural
class of metric \Ss s \cite{FHK}, \cite{Ha1}, \cite{He}. The discussion of
higher order \Ss s and higher order mean \df\ \q s
on regular subsets of Euclidean spaces is also possible \cite{16'} in
the context of the generalized Taylor--Whitney jets,
\cite{Boj13}--\cite{BIK}, \cite{BBM2}, \cite{Gla1}, \cite{Wh3}.
This paper is a direct sequel to \cite{r0}, \cite{Boj-gr}.
\end{abstract}

\section*{Introduction}
Averaging and integration, especially integration by parts formulae and their
multidimensional analogues (Stokes--Green theorems etc), lie at the heart of
Sobolev's theory and his concepts of duality in linear function spaces.

The Hardy--Littlewood idea \cite{HL} of maximal function gives the natural
and most direct way to introduce and characterize the Sobolev spaces
$W^{m,p}(\R^{n})$ and $W_{\rm loc}^{m,p}(\R^{n})$ and their far going
generalizations to measure metric spaces \cite{Ha1}, \cite{He}, \cite{HeK2}
etc.

We illustrate this fact and the corresponding ideas on the  simplest case of Sobolev spaces
$W^{m,p}(\R)$, $m\ge 1$, $p\ge 1$ ($m$--an integer) on the real line, starting with the case
$m=1$, $p=1$. Since a function
$f\in W^{1,1}(\R)$ is absolutely continuous,  by the fundamental theorem of the calculus for
$x,y$ real, $x<y$, we have
\beq\label{0.1}
f(y)-f(x) = \int_{x}^{y}f'(t)\,dt.
\e
Hence for $x<\xi<y$
\beq\label{0.2}
f(y)-f(x) = (y-x) \left[\frac{\xi-x}{y-x}\kint{x}{\xi}f'(t)\,dt + \frac{y-\xi}{y-x}\kint{\xi}{y}
f'(t)\,dt\right]
\e
and, by the definition\footnote{We skip here over the precise definitions, distinctions of left
$M_{L}$ and right $M_{R}$
maximal functions on the line etc., for details see \cite{BKZ}, \cite{Gra}, \cite{Ha1}.} of the Hardy--Littlewood maximal
function \cite{HL}, \cite{Gra}
\beq\label{0.3}
g_{f}(x) \equiv M(|f'|)(x)~\footnote{Here, and below:
$\fint{V}{} g\,dx \equiv \frac{1}{|V|}\int_{V}g\,dx$, $|V| > 0$.} = \sup_{\xi>x} \kint{x}{\xi} |f'(x)|\,dt
\e
we get the pointwise inequality
\beq\label{0.4}
|f(y)-f(x)|\le |y-x| [g_{f}(x) + g_{f}(y)].
\e
Now for $f\in W_{\rm loc}^{1,p}(\R)$, $p>1$, $f'(t)\equiv \nabla f(t)$, \cite{HL}, \cite{Gra}, \cite{St1}
\beq\label{0.5}
\|g_{f}\|_{p}\le C_{p} \|\nabla f\|_{p}
\e
and we conclude \cite{Boj13}, \cite{BH}, \cite{Ha1} that, for $p>1$, the pointwise inequality \eqref{0.4}
characterizes the Sobolev class $W^{1,p}(\R^{n})$ \cite{BH}, \cite{Ha1}, \cite{Ha2}.

It is a classical fact of analysis that, for $f\in W^{1,p}(\R)$, $p>1$, \eqref{0.1} implies, by H\"older inequality, the following estimate
\beq\label{0.6}
|f(y)-f(x)|\le |y-x|^\alpha \|f'\|_{p}\quad\ \alpha=1-\frac1p > 0.
\e
This describes, when generalized to $\R^{n}$, $n>1$, the simplest of the famous Sobolev imbeddings
\beq\label{0.7}
 W_{\rm loc}^{1,p}(\R^{n})\to C_{\rm loc}^{\alpha}(\R^{n}), \quad\ \alpha=1-\frac{n}{p},\ p>n,\ n\ge 1.
\e
For $p\le n$ the Sobolev functions  admit discontinuities at, possibly everywhere
dense, sets of points and may be even locally unbounded \cite{BIN}.

The inequalities \eqref{0.6} and \eqref{0.7} imply the uniform continuity of bounded families of functions
in $W^{1,p}(\R^{n})$, $p>n$, and play a crucial role for various compactness arguments in the calculus of variations,
PDE, and many other applications of Sobolev spaces in analysis and geometry.

However \eqref{0.6} is very far away from characterizing the Sobolev functions.
In particular the  Weierstrass
nowhere differentiable function \cite{EG}, \cite{Gra}, \cite{FHK} is H\"older continuous, for some $\alpha>0$.
In fact a much deeper statement holds: in general an $\alpha$-H\"older continuous function, with
$\alpha<1$, cannot
be ``corrected'' on any subset of positive measure of its domain of definition to become
 Lipschitz (smooth)!

On the contrary,  any function $f$ in the Sobolev class $W^{1,p}(\R^{n})$, or any function just satisfying the inequality \eqref{0.4}, can be modified on a subset of arbitrary small positive
Lebesgue measure in its domain of definition
 to become Lipschitz or even $C^{1}$ \cite{EG}. This is the important Luzin's property (or Luzin approximation property)
\cite{Boj-f}, \cite{EG}, \cite{HdN}, \cite{BIK}, of function classes, satisfying inequalities of the type \eqref{0.4} and their numerous generalizations.

The essence of the argument runs as follows: For $L>0$ consider the
set $E_{L}=\{x: g_{f}(x)\le L\}$.
On $E_{L}$ the function $f$ is (uniformly) Lipschitz. By Tschebyscheff inequality the Lebesgue measure of the complement $CE_{L}$ of $E_{L}$
\beq\label{0.8}
|CE_{L}| \to 0\quad\ \text{for}\ L\to \infty
\e
Invoking the Tietze--McShane--Kirschbraun--Whitney extension theorems \cite{EG},
\cite{HdN}, \cite{He} we come to the
required conclusions.

For some Sobolev type function spaces the details of this procedure are given in the
quoted references. Some others will be mentioned in Section 3 below.

Let us remark at this point that the remarkable Luzin's property has
essentially a semi global character: it refers to the behaviour of the
function (mapping!) on the \emph{complement} of a set \emph{very small in
measure}!  The deep structure theorems of N. N. Luzin for measurable
functions on general metric spaces allow us to understand the duality between
measurability and continuity: Luzin's continuity of a mapping is an
infinitesimal expression of the global property of measurability of the
mapping!!  In a still more refined form these questions led A.~Denjoy and
A.~Khintchine,
\cite{Federer}, to the discovery of the concept of the approximate continuity
almost everywhere as Luzin's dual to measurability. See e.g.\ references
in~\cite{Boj-f}.
This point of view is presented in extenso in \cite{Boj-f}.

When applied to the Sobolev type function spaces it allowed Hassler Whitney,
\cite{Wh2}, \cite{Wh4'}, to extend the Luzin's duality to continuously
differentiable functions in the class $C^{1}(\R)$, on the real line. In
\cite{Wh4'} H. Whitney showed that the semiglobal Luzin's property of
order~$1$, see \cite{Boj-f}, has an infinitesimal description in terms of the
almost everywhere approximate Peano differentiability \cite{Federer}. In this
connection let us also state that the long
standing Federer--Whitney conjecture, (3.1.17) on pp.~228--229 of
\cite{Federer}, was solved in \cite{Boj-f} thanks to the work of F. Ch. Liu
and the author described in \cite{Boj-f}\,\footnote{At this point it is
necessary to say that the Russian edition translator's remark on p.~247 of
the Federer's monograph \cite{Federer} (Nauka, Moscow, 1987) and the referred
theorem of E. E. Movskovitz is erroneous.}. These works allowed us also to
understand the case of Luzin's property for Sobolev spaces $W^{m,p}(\R^{n})$
of higher order $m\ge 1$. The detailed discussion is postponed to \cite{16'}.

Another important property of the pointwise inequalities
\eqref{0.4}--\eqref{e9} (below, \eqref{e9} for
the euclidean spaces $\R^{n}$)
is that they are directly related, in fact equivalent, with the classical Riesz potential estimates:
the starting point of  Sobolev's theory \cite{Sob}, \cite{Sob3}, \cite{Sob2} in $\R^{n}$ is,
essentially,  again the fundamental theorem of the calculus \eqref{0.1} written in the form
\beq\label{0.9}
f(x)-f(y) = -\int_{0}^{|x-y|} D_{t}f(x+t\omega)\,dt
\e
where $x, y$ in $\R^{n}$ are in the ball $B_{r}\equiv B(x,r)=\{y, |x-y|<r\}$
and $\omega=\frac{y-x}{|y-x|}$ is the unit vector
in $\R^{n}$. Averaging \eqref{0.9}  with respect to $y$ over the ball $B_{r}$,
after introducing polar coordinates and a change of
variables, Sobolev comes to the pointwise estimates
\beq\label{0.10}
|f(x)-f_{B}|\le C(n)I_{1}(|\nabla f|)(x)
\e
where $I_{\alpha}(g)(x)\equiv (|y|^{\alpha-n}* g)(x)$ is the Riesz potential of order $\alpha>0$. In \eqref{0.10}
$\alpha=1$.

Written in a symmetric way \eqref{0.10} takes the form
\beq\label{0.11}
|f(x)-f(y)| \le C(n)\left(I_{1}(|\nabla f|)(x) + I_{1}(|\nabla f|)(y)\right),
\e
which is a special case of \eqref{0.4} in the classical Euclidean case of $\R^{n}$.

As a matter of fact, the inequality
\eqref{0.10} was the starting point of the systematic approach to pointwise Sobolev inequalities initiated in \cite{BH}.
For the inverse way: from \eqref{0.4} to the pointwise Poincar\'e type inequality
\eqref{0.10} see \cite{Ha1}, \cite{Ha2},
\cite{HaK2}, \cite{Kra} and the inequalities \eqref{e11}, \eqref{e12} below.

S. L. Sobolev introduced his spaces $W^{m,p}(\R^{n})$, $m\ge 1$, $p\ge 1$
($W^{m,p}(\Omega)$, $\Omega$ a domain in $\R^{n}$)
as closed subspaces of the Lebesgue spaces $L^{p}(\R^{n})$, $(L^{p}(\Omega))$ defined by some integral identities resulting from the
classical integration by parts formulae \cite{Sob1}, \cite{Sob}, \cite{Sob3}, \cite{BIN}, \cite{FHK}.

These identities also imply the recurrence relations (under some obvious conditions on~$f$)
\beq\label{0.12}
f\in W^{m,p}(\R^{n})\iff \nabla f\in W^{m-1,p}(\R^{n}),\quad
(m\ge 1),\ (W^{0,p}(\R^{n})\equiv L^{p}(\R^{n}))
\e
where $\nabla f$ is the generalized (weak) Sobolev
gradient of~$f$. They open the way for the inductive, with respect to $m$,
 treatment of the Sobolev spaces of order $m\ge 2$~\cite{88}.

However, for various reasons, the direct approach to $W^{m,p}(\R^{n})$ (for $m\ge 2$) has been, and continues to be,
widely used as well \cite{BIN}, \cite{St2}, \cite{Tri}.
The gradient operators $\nabla^{k}$ $(k\ge 1)$ combined with the finite difference and shift operators
\beq\label{0.12a}
\Delta f = \Delta_{y}f = f(y)-f(x),\quad\ T_{\tau}f(x)\equiv f_{\tau}(x)=f(x+\tau),\ y,x,\tau\in \R^{n}
\e
and their iterates create the general framework for the discussion of Sobolev type functions spaces \cite{BIN},
\cite{r0}, \cite{Sob2}, \cite{St1} in general.

In the immense and permanently growing literature of the subject there are thousands of works describing Sobolev spaces
as Banach space closures of smooth functions in the corresponding Sobolev integral norms.

In various constructive computational treatments of PDE and their applications in sciences many new, ingenious
methods were invented \cite{DS}, \cite{Shv2}, \cite{St1} etc.

Sobolev functions on the real line $\R^{1}$ occupy a
special place in the general theory of  Sobolev spaces $W^{m,p}(\R^{n})$:
elements $f$ in $W^{m,p}(\R^{n})$ (or on manifolds)
have representatives $\overline f$, $f\equiv \overline f$ a.e.,
\cite{88}, \cite{BIN}, which, restricted to almost all lines $l\subset \R^{n}$ parallel to the coordinate axes, may be regarded as generalizations of Sobolev functions in $W^{m,p}(\R^{1})$, $f_{l}=\overline f_{|l}$ is in $W^{m,p}(l)$.

For a short review of the theory of Sobolev spaces on the real line, with
special emphasis on the use of the concepts
of Hardy--Littlewood maximal functions, pointwise inequalities, $B$-Splines, numerical analysis and
approximation theory \cite{DS}, \cite{EG}, \cite{Fi}, \cite{Shv1}, see \cite{BKZ} where the presentation is somehow
adapted to the needs and concepts of this paper (see also \cite{Boj-gr}).

For $m=1$, $p=1$ the space $W^{1,1}(\Omega)$, $\Omega\in \R^{n}$ is
identified with the famous class $ACL(\Omega)$. The class ACL is widely used
in the analytic theory of quasiconformal mappings, e.g.\ \cite{BIw},
\cite{He}, \cite{He2}, of real valued functions absolutely continuous on
almost all lines (line segments) in $\Omega$, parallel to the coordinate axes,
and whose partial derivatives $f_{x_{i}}$, defined a.e. in $\Omega$, belong
to $L^{p}(\Omega)$.

This remarkable property, the directional regularity of Sobolev functions is
analogous to the famous Hartogs theorem, see \cite{JP}, (and references
therein) for separately holomorphic functions in $\C^{n}$.  It is a kind of
Fubini theorem for Sobolev functions. The general case of the theorem, for
the product decomposition $\R^{n}=\R^{k}\times \R^{n-k}$, $k\ge 1$ and for
manifolds, is due to S.~M.~Nikolskii in his  early paper \cite{Nik}, see also
\cite{BIN}, \cite{Boj4}, \cite{88}.

The pointwise inequalities characterizing the elements $f\in W_{\rm loc}^{m,p}(\R^{n})$,
$m\ge 2$, are formulated, for a
pair $(x,y)$, $x\ne y$, of points in $\R^{n}$ on an affine line $l=l(x,y)$ in $\R^{n}$, as estimates of the type
\eqref{0.4} for the ``errors'' or differences
\beq\label{0.13}
R^{m-1}f(x,y) \equiv R_{Z}^{m-1}f(x,y)=f(y) - L(y;f,x_{0},\ldots,x_{m-1}),\quad x_{0}=x
\e
between the value $f(y)$ and the value at $t=y$ of the Newton--Lagrange interpolation polynomial
$p_{n}(t) = L(t;f,x_{0},\ldots,x_{m-1})$ \cite{Boj-gr}, \cite{r5}, \cite{DS}, \cite{Ge}, \cite{Shv1} of order $n=m-1$,
interpolating the function $f(t)$ at the nodes $x_{i}$
\beq\label{0.14}
p_{n}(x_{i}) = f(x_{i})\quad\ i=0,\ldots,n.
\e
The nodes $x_{i}$ in \eqref{0.13} and \eqref{0.14} are assumed to belong to the affine line $l=l(x,y)$,
or rather to the affine segment $I(x,y)\subset l$; $t$ is the natural parameter on~$l$.

The letter $Z$
in \eqref{0.13} stands for the general term of ``interpolation scheme'' \cite{4'}, \cite{Gra} i.e.\enskip the set of general
conditions of the type of Hermite--Birkhoff--La\-grange conditions: see \cite{4'}, (or their equivalence classes under
affine transformations of $\R^{n}$) used to define the interpolating polynomial $p_{n}(t)$.

Though the classical interpolation theory \cite{4'}, \cite{r5}, \cite{DS}, \cite{EG} is considered for complex
valued functions on the real line $\R^{1}$ only, the term Newton--Lagrange interpolation polynomial $L(t;f,x_{0},\ldots,x_{m-1})$
is meaningful for affine lines in~$\R^{n}$. See the comments in this respect in \cite{r0}, \cite{Boj-gr}.

For simplicity, in this paper, we assume that all nodes $x_{i}$ are simple and ordered
\beq\label{0.15}
x=x_{0}<x_{1}<\cdots < x_{n}\le y
\e
except the case $x=x_{0}=\cdots =x_{n}<y$ when the Lagrange interpolation
 reduces to the Taylor formula.

For many reasons it is convenient to represent the interpolating polynomial
in the Newton form \cite{Boo}, \cite{22'}, \cite{St1}
\beq\label{0.16}
p_{n}(t) = \sum_{i=0}^{n} c_{i}q_{i}(t),\quad\ q_{0}(t)\equiv 1
\e
where the polynomials in $t$, $q_{i}(t) = \sum_{k=0}^{i=1} (t-x_{k})$, are the Newton polynomials for the assigned nodes $x_{i}$, $i=0,\ldots,
n-1$ of the interpolation scheme~$Z$. The coefficients $c_{i}$ in \eqref{0.16} are the divided differences
of the function $f$ for the scheme~$Z$. In the classical notation \cite{4'}, \cite{r7}, \cite{r8}, \cite{Schu},
\beq\label{0.17}
c_{i}=f[x_{0},\ldots,x_{i}]
\e
and the remainder $R_{Z}^{n}f(t,x)$ takes the form
\beq\label{0.18}
R_{Z}^{n}f(t,x) = f[x_{0},\ldots,x_{n},t]q_{n+1}(t)\quad\ x_{0}=x.
\e
For Sobolev spaces $W^{m,1}(\R^{1})$, $m=n+1$, the normalized remainders have the beautiful form of
the general formula
\beq\label{0.19}
\frac{R_{Z}^{m-1}f(t,x)}{q_{m}(t)} = f[x_{0},\ldots,x_{n},t]
\e
i.e. they are represented as divided differences themselves.

For suitable choices of the Lagrange--Newton interpolation scheme,
the formula \eqref{0.18} comprises all
classical formulas for the Newton--Taylor--Lagrange remainders used in the literature \cite{Schu}.
Thus, for $m=1$ and $f\in C^{1}(R)$, $t=y$, we get
the identity
\beq\label{0.20}
Rf(y,x) = R^{0}f(y,x) = (y-x)f[y,x] = (y-x)\int_{0}^{1} f'(x+ \tau(y-x))\,d\tau
\e
used in \eqref{0.2} and \eqref{0.4}.

For $m\ge 2$ and the equidistant nodes $x_{i}=x+ih$, $h=\frac{y-x}{m}$, $n=m-1$, we get the
classical finite
difference remainders \cite{r0}, \cite{Boj-gr}, \cite{Tri}
\beq\label{0.21}
R^{m-1}f(y,x) = \Delta^{m} f(y,x) \equiv \sum_{i=0}^{m}(-1)^{m-i} \binom{m}{i}f(x_{i})
\e
leading to the pointwise inequalities
\beq\label{0.22}
|\Delta^{m} f(y,x)| \le |y-x|^{m}[g_{f}(x) + g_{f}(y)]
\e
characterizing the Sobolev spaces $W^{m,p}(\R^{n})$, $p>1$ by the condition $g_{f}\in L^{p}(\R^{n})$. Thus \eqref{0.22}
is a direct generalization of \eqref{0.4} and \eqref{e9} below.

If the nodes $x_{i}$, $i=1,\ldots,n$ shrink to the center $x_{0}=x$,
\beq\label{0.23}
\lim_{l\to \infty}x_{i}^{l}\to x\quad i=1,\ldots,m \ \text{or}\ \|y-x\|\to 0
\e
for a family $Z^{l}$ of interpolation schemes, the remainders $R_{Z^{l}}^{m-1}f(y;x)$ reduce to the classical Taylor or, Taylor--Whitney  \cite{Wh2} remainders for $f$ differentiable (of class $C^{m}(\R^{n}))$. We recall that for $f\in C^{m}(\R^{n})$ $(x_{m}=y)$
\beq\label{0.24}
f[x,x_{1},\ldots,x_{m}]\to \frac{f^{(m)}(x)}{m!}
\e
when $|x_{i}-x|\to 0$, $i=1,\ldots,m$ \cite{r7}. In this case
 the Newton representation formula \eqref{0.16} reduces to the classical Taylor formula \cite{BKZ}, \cite{r7}.

The general form of the pointwise inequality for functions $f$ in $W^{m,p}(\R^{m})$ has the form
\beq\label{0.25}
|f[x,x_{1},\ldots,x_{m-1},y]|\le g_{f}(x)+ g_{f}(y),
\e
for some $g_{f}\in L^{p}(\R^{n})$.

The function $g_{f}$ is majorized from above by the Hardy--Littlewood maximal function of the $m$-th Sobolev
gradient $\nabla^{m}f$. For $m=1$, it is shown below that $g_{f}$ can be taken as $g_{f}=MQf$, the mean maximal quotient,
see \eqref{e9}. For the (harmonic) interpolation scheme with equidistant nodes \eqref{0.22} the proof of the asserted
estimate is given in \cite{Boj-gr}.

The general case, for functions on the real-line \cite{r8}, \cite{r7},
\cite{Schu}, is deduced from Genocchi--Hermite formula, \cite{BKZ},
\cite{Schu}, in analogy with the direct proof of \eqref{0.4} above. For
Taylor--Whitney interpolation scheme \eqref{0.23}--\eqref{0.24} the
inequality \eqref{0.25} reduces to the case considered in \cite{BH}.

The discussion of the simplest, now available, proof of \eqref{0.25} for the
general interpolation schemes with arbitrary nodes $\{x_{i}\}$ will be given
in \cite{16'}.

Since the Newton polynomial $q_{m}(y)$ for the interpolation
scheme $Z: \{x_{i}
\le y,\, i=1,\ldots,m-1, x_{0}=x\}$ is majorized by $|y-x|^{m}$ we can also
write \eqref{0.25} in the more familiar form
\beq\label{0.26}
|R^{m-1}f(y,x)| \equiv |R_{Z}^{m-1}f(y,x)| \le |y-x|^{m} (g_{f}(x)+ g_{f}(y)).
\e

For $m=1$ \eqref{0.26} reduces to the inequality
\beq\label{0.27}
|Rf(y,x)| = |f(y) - f(x)|\le |y-x| [g_{f}(x)+ g_{f}(y)]
\e
which also makes sense for measure-metric spaces
$(X,d,\mu)$, $x,y\in X$ with distance $d(x,y)\doteq |y-x|$ and measure~$\mu$.

Of course the coefficients $g_{f}$ in the formulas \eqref{0.22}--\eqref{0.27}
vary from case to case and they appear in the literature on Sobolev spaces
and their applications, over the last 25 years, at various places with various names used: variable
Lipschitz coefficient, Sobolev metric gradient etc. Perhaps an appropriate term for $g_{f}$ could be
\emph{Sobolev smoothness} or \emph{Sobolev smoothness density}, which occasionally is used in this and
related papers, of the author and his students and colleagues \cite{Boj13},
\cite{BH}, \cite{BHS}, \cite{Ha1}, \cite{Ha2}, \cite{HdN}, \cite{HeK2}, etc.

In the works involved, the right hand side coefficients $g_{f}(x)$, in general, are
estimated from above in terms of the maximal functions of the $m$-th gradient $|\nabla^{m}f|$
evaluated at $x$ and some geometric constants (universal, depending on $n$,
$p$ and $Z$ only) in analogy with \eqref{e8} below.

Even for functions on the real line the inequalities \eqref{0.22} and \eqref{0.25}
seem to be absent in the literature.

This last point has been confronted with the expertise at some leading mathematical
centers in the area (Sobolev Institute at Novosibirsk, Courant Institute at NYU,
some other universities in USA
and Europe, especially Finland, and a number of leading experts in Sobolev space theory.

The author will be grateful for any additional related references.

A natural and important problem is the understanding of the
dependence of the pointwise estimate \eqref{0.25} on the interpolation
scheme $Z$ \cite{4'}, \cite{Gla1}. It
turns out that the divided difference $f[Z]$ \eqref{0.25}, and \eqref{0.22}, \eqref{0.26} as well,
 are estimated from above by the values of the Sobolev smoothness $g_{f}$ at the extremal points $x$ and $y$ of $Z$, and
the upper bound does not depend on the intermediate nodes $x_{i}$ in \eqref{0.15}.

Even in the  intensively
studied classical case of the strictly one dimensional interpolation and approximation theory on the
real line $\R$ \cite{r5}, \cite{DS}, \cite{EG}, \cite{GL}, \cite{Schu}\dots\ this problem seems
 not to have been
satisfactorily clarified. Also in the rather recent papers of
H. Triebel \cite{Tri} and his school \cite{HT}
(and some others following, available on the net) the independence
of the right hand side in \eqref{0.22} from the intermediate nodes is
overlooked.

Notice, in this connection, that in  \cite{DS} a somehow different concept of A. P. Calderon's maximal function
is used to characterize the higher order Sobolev type smoothness of measurable functions. Here again the values at
all intermediate nodes of the corresponding maximal functions are used in the estimates of the expressions \eqref{0.21}.
Apparently, in view of \eqref{0.22} and \cite{Boj-gr}, this fact deserves to be reconsidered.

The pointwise inequalities
\eqref{0.22} and \eqref{0.27} are intimately related with S.~M.~Nikolskii's
trace theorem quoted above:
all quantities in the left hand sides of \eqref{0.24}, \eqref{0.25}, \eqref{0.26} are
determined by the fact that the restriction $f_{I}$ of the function $f$ to the segment $I=I(x,y)$
belongs to the Sobolev space $W^{m,p}[I]$ on the segment~$I$.
The evaluation of the terms $g_{f}(x)$ and $g_{f}(y)$ in the right hand sides
of \eqref{0.22}--\eqref{0.26} is
determined by averaging of the gradients, difference quotients, and other quantities
describing the ``generalized'' Sobolev regularity of $f$,
over open neighbourhoods of the segment $I(x,y)$ in the enhancing space $\R^{n}$ see
\eqref{e2} and \eqref{e16}, \eqref{eq2.7} below.
The quantities $g_{f}(x)$ and $g_{f}(y)$ reflect the fact that the function $f\in W^{m,p}(\R^{n})$.

All in all, the pointwise inequalities appear as primary facts of the Sobolev space theory:

\smallskip
they are formulated and proved at the very initial stages of the theory;

in a properly understood conceptual context they have simple and natural
proofs~\footnote{in mathematics simple ideas usually come last (J. Hadamard).};

in a rather direct way they lead to the central results of the theory:
the famous Sobolev embedding theorems and delicate, far from obvious, continuity,
differentiability, Peano approximate differentiability, Luzin's property
 and other crucial properties of Sobolev space theory
and allow us to see them in their natural general settings.
\smallskip

For the Sobolev space $W^{m,p}$ these are: for $m=1$ measure, metric spaces \cite{Ha1} and for $m\ge 2$ measurable subsets or submanifolds of
domains in $\R^{n}$.

The text below and \cite{16'} is a continuation of a series of papers published over
the last two decades, partially quoted in our references,
illustrating some of the above statements.

\section{Maximal mean difference quotients}
Let $f$ be a real valued locally integrable \f\ $f\in L\loc(\R^n)$. The
maximal mean \df\ \q\ at the point $x\in\R^n$, denoted by $MQf(x)$, is
defined as the following l.u.b.
\beq\label{e1}
MQf(x)=\sup_{r>0} M_rQf(x)
\e
of the \a ed values of the \df\ \q s
\beq\label{e2}
M_rQf(x)=\fint{B(x,r)}{} \frac{|f(z)-f(x)|}{|z-x|}\,dz.
\e
(Here $B(x,r)$ is the open ball $\{y: |y-x|<r\}$.)

Our general assumption is that for bounded subsets $\O\subset \R^n$ the
integral
\beq\label{e3}
\int_\O \int_\O \frac{|f(z)-f(x)|}{|z-x|}\,dz\,dx<\iy
\e
is finite.
Condition \er3 implies that the \a ed values \er2 are well defined for $x$
a.e.

\bdf \label{d1}
For $p\ge1$ the class $MQ^{1,p}(\R^n)$ is defined as the set of all $f$ in
$L^p(\R^n)$ for which $MQf$ is in $L^p(\R^n)$.
\edf

If $MQf(x)$ is finite a.e.\ the following proposition holds

\bpr \label{p1}
With an absolute \ct\ $c=c(n)$ the pointwise \iq y
\beq\label{e4}
|f(x)-f(y)|\le c|x-y|(MQf(x)+MQf(y))
\e
 holds a.e.; specifically
for all $x,y$ for which all terms in \er4 are meaningful and finite, e.g.\
for all Lebesgue points of~$f$ $($the right-hand side of \er4 is defined,
possibly $\iy$, for all $x,y,\ x\ne y)$.
\epr

\proof
For $x,y$ fixed let $\S_r=\S_r(x,y)$ be the spherical segment
\beq\label{e5}
\S_r=B(x,r)\cap B(y,r), \quad r=|x-y|.
\e
For all $z\in\S_r\ z\ne x\ne y$ we have
\beq\label{e6}
|f(x)-f(y)|\le |f(x)-f(z)|\,\frac{|z-x|}{|z-x|} +|f(z)-f(y)|\,\frac{|z-y|}{|z-y|}
\e
thus
\beq\label{e7}
\frac{|f(x)-f(y)|}{|x-y|}\le \frac{|f(x)-f(z)|}{|x-z|}
+\frac{|f(z)-f(y)|}{|z-y|},
\e
since $|x-z|\le|x-y|$, $|y-z|\le|x-y|$.

\E\a ing \er7 \wrt $z$ over $\S_r$ we get
\begin{multline*}
\frac{|f(x)-f(y)|}{|x-y|}\\ \le \frac{|B(x,r)|}{|\S_r|}\biggl(
\fint{B(x,r)}{}\frac{|f(x)-f(z)|}{|x-z|}\,dz
+\fint{B(y,r)}{}\frac{|f(z)-f(y)|}{|z-y|}\,dz\biggr).
\end{multline*}
By the geometry of $\R^n$
\beq\label{e8}
\frac{|B(x,r)|}{|\S_r|}=\frac{|B(y,r)|}{|\S_r|}=c(n)
\e
and \er4 follows with $c=c(n)$.

Written in the form
\beq\label{e9}
|f(x)-f(y)|\le |x-y|(g(x)+g(y)), \quad g(x)=g_f(x)=MQf(x),
\e
\er4 is the first and, apparently, most natural example of the pointwise \iq ies
characterizing the general classes of \f\ spaces of Sobolev type, see \cite{BH},
\cite{Ha1}.

We shall call the \iq y \er9 the \ti{Sobolev \pw\ \iq y} of order~1.

\E\iq ies \er9 (and some their modifications) have been widely used in the
last two decades to define various generalizations of classical \Ss s
$W^{1,p}(\R^n)$ \cite{r0}, \cite{Ha1}, \cite{HaK1}, \cite{Kos}, \cite{He}.
We shall use for these spaces the general
term \emph{metric Sobolev spaces} \cite{He}, \cite{HeK2} and, for a given $f\in L^p(\R^n)$, any admissible
\f~$g$ in the \rh\ of~\er9 will be termed a Sobolev metric gradient of~$f$.
The collection $SMG(f)$ of all \Sm \g s for $f\in L^p(\R^n)$ is a convex
(closed) subset of nonnegative measurable \f s in~$\R^n$. In fact it is a
lattice (i.e., $g_1$ and $g_2\in SMG(f)$ implies $\min(g_1,g_2)\in SMG(f)$).
Also if $g\in SMG(f)$ any $g_1\ge g$ is in $SMG(f)$ as well.

In particular, the following definition \cite{Ha1} is broadly used.

\bdf \label{d2}
$M^{1,p}(\R^n)$ is the class of all $f\in L^p\loc(\R^n)$ which admit \Sm \g s
in $L^p(\R^n)$.
\edf

Our Proposition \ref{p1} implies

\bco \label{c1}
 $MQ^{1,p}(\R^n)$ $(p\ge1)$ is a \Sm space: $MQ^{1,p} \subset M^{1,p}(\R^n)$.
\eco

Notice that each $f\in MQ^{1,p}$ has a well specified \Sm \g\ $g=MQ(f)$
defined by the formulas \er1 and \er2.

In fact the above definitions and the Hardy--Littlewood theorem (recalled
below) imply also

\bpr \label{p2}
If $f\in M^{1,p}(\R^n)$, $p\ge 1$, then any $g \in SMG(f)$ can be estimated from
below pointwise by the mean maximal quotient
\beq \label{e9.0}
MQf(x)\le g(x)+Mg(x)\le 2Mg(x)\quad \hbox{for a.e. }x,
\e
where $Mg(x)$ is the \HL\ maximal \f\ for $g\in L^p\loc(\R^n)$, $g\ge 0$.
\epr

This is a direct consequence of the definitions and the \iq y \er9.

In particular, the next corollary follows.

\bco \label{c2}
For $p>1$, $M^{1,p}(\R^n)\subset MQ^{1,p}(\R^n)$ $($with equivalent
semi\-norms$)$.
\eco

In the sense described by Proposition \ref{p2}, the maximal mean quotient
takes the role of a ``\pw\ minimal'' element in $SMG(f)$.

$MQ^{1,p} (\R^n)$ is a linear space with the seminorm $\sn f_{1,p}
=\|MQf\|_{L_p}$.

In Proposition \ref{p1} the Sobolev pointwise \iq y \er9 was obtained by
averaging the \df\ \q\ \er7. If, instead, we \a e the oscillation
$|f(z)-f(x)|$  on open non empty convex subsets $\S$ of~$\R^n$ we come to the \P\ \iq
ies as we now recall.

For the pair of points $x,z$ of $\S,\ x\ne z$ we have
\beq \label{e9.1}
|f(z)-f(x)|\le |z-x|\,\frac{|f(z)-f(x)|}{|z-x|}
\le\diam \S\,\frac{|f(z)-f(x)|}{|z-x|}.
\e
If $\S$, rather a family $\S (r)$, $\S = \S(r)$, is a subset family of $B(x,r)$ such that
\beq \label{e9.2}
\frac{\,|B|\,}{|\S|}\le C, \quad \hbox{independent of }r,
\e
e.g., $\S=B(x,r)$ or $\S=B(x,r)\cap B(y,r)$, we get
\[
|f_\S -f(x)|\le \diam \S\,\frac{\,|B|\,}{|\S|} \fint B {}\frac{|f(z)-f(x)|}
{|z-x|}\,dz\le C(\diam\S) M_rQf(x)
\]
consequently
\beq \label{e9.3}
|f_\S -f(x)|\le C(\diam\S) MQf(x).
\e
In particular, for $\S=B(x,r)$
\beq \label{e9.4}
|f(x)-f_{B(x,r)}|\le Cr MQf(x) \stackrel{.}{\equiv} rg(x)
\e
for almost all $x$, e.g.\ for all Lebesgue points of $f\in L^p(\R^n)$ and,
analogously, for $\S_r=B(x,r)\cap B(y,r)$, $c_n$~as in \er8,
\beq \label{e9.5}
\bigl|f(x)-f_{\S_r}\bigr| \le rc_ng(x), \quad r=|x-y|.
\e
\er{9.4} and \er{9.5} can be viewed as the pointwise form of the \P\ type \iq
ies for the \f\ space $MQ^{1,p}$.

Raised to the exponent $p\ge1$ and integrated, \er{9.4} and \er{9.5} lead to
the familiar integral form of the \P\ \iq ies. They are recalled below for the
case of smooth \f s $f\in C^1(\R^n)$ or $f\in W^{1,p}(\R^n)$.

For \f s in the class $C^1(\R^n)$ the fundamental theorem of calculus
\beq\label{e10}
f(z)-f(x)=\int_0^1 \frac d{dt}f(x+t(z-x))\,dt=\int_0^1
\langle \n f(x+t(z-x)),z-x\rangle\, dt
\e
implies the \iq y
\beq\label{e11}
|f(z)-f(x)|\le |z-x|\int_0^1 |\n f|(x+t(z-x))\,dt.
\e

Let $\S$ be a convex set (in $\R^n$) and consider pairs $x,z\in\S$, then
$|x-z|\le \diam\S$. \E\a ing \er{11} over~$\S$ \wrt $x$ and~$z$ we get
\begin{align}
|f(x)-f_\S|&\le \diam\S \int_0^1 \fint \S{}{} |\n f|(x+t(z-x))\,dt\,dz\notag\\
&\le c_n\diam\S \fint \S{}{} |\n f|\,dz\label{e12}
\end{align}
\beq
\fint \S{} |f(x)-f_\S|\,dx \le \diam\S \int_0^1 \fint \S{}{}
|\n f|(x+t(z-x))\,dt\,dz\,dx. \label{e13}
\e
After standard computation using the change of variables $\z=x+t(z-x)$, $d\z
=t^n\,dz$ and the H\"older \iq y $\fint \S{}\ g\,d\mu\le \bigl(\fint \S{}\
g^p\,d\mu\bigr)^{1/p}$, $p\ge1$, $g\in L^p(\S,\mu)$, we come \cite{GT}, \cite{Kos}
to the classical
Poincar\'e \iq ies
\beq\label{e14}
\fint \S{} |f-f_\S|\,dx \le c_{n,p}\diam \S \fint \S{}|\n f|\,dx
\e
and
\beq\label{e15}
\biggl(\fint \S{} |f-f_\S|^p\,dx\biggr)^{1/p} \le c_{n,p}\diam \S
\biggl(\fint \S{}|\n f|^p\,dx\biggr)^{1/p}, \quad p\ge 1.
\e

The \P\ integral \iq ies \er{15}, \er{14}, \er{13} are the only ones
above\footnote{\f s in $W^{1,p}$ are discontinuous and have values a.e. only}
which are meaningful and valid for arbitrary \f s~$f$ in the
classes $W^{1,p}\loc(\R^n)$, $p\ge1$. They are generally considered to be
much weaker than the more subtle Sobolev integral \iq ies discussed later.
Sobolev integral
\iq ies have an essentially different form for the values of the parameter~$p$:
$1\le p<n$, $p=n$ and $p>n$.

The \pw\ form \er{9.3}, \er{9.4} of the \P\ \iq
ies can be given sense for a.e.\ point of~$\S$ only (for all values of the parameter
$p\ge 1$) if the gradient $|\n f|$
in the \rh\ of the \P\ \iq ies in $W^{1,p}(\R^n)$ is replaced by the more
refined concept of the Hardy--Littlewood maximal
\f\ of the gradient $|\n f|$, as will be explained later.

Actually Sobolev was \a ing the gradients, not the \df\ \q s, but his
fundamental papers \cite{Sob1}, \cite{Sob}, \cite{Sob3} led to the discovery of concepts, ideas and
new mathematical facts of unprecedented importance for higher analysis,
differential and partial differential \eq s and applications.

In $\R^n$ the relation of the class $MQ^{1,1}(\R^n)$ with the
class of smooth functions and classical \Ss s $W^{1,1}(\R^n)$ is described by
the elementary (well known)

\blm \label{l1}
For $f\in C^1(\R^n)$ and any ball $B(x,r)$ the \iq y holds
\bml \label{e16}
\fint{B(x,r)}{}\frac{|f(z)-f(x)|}{|z-x|}\,dz=\int_0^1
\fint{B(x,r)}{}|\nabla f(x+(z-x)t)|\,dz\,dt\\{}\le
\int_0^1 \fint{B(x,tr)}{}|\nabla f(z)|\,dz\,dt.
\e
\elm

\proof
Use \er{11} and the change of variable: $\z=x+(z-x)t$, $d\z=t^n\,dz$,
transforming the ball $B(x,r)$ into $B(x,tr)$.
\endproof

For \f s in the classical \Ss\ $W^{1,p}\loc(\R^n)$, $1\le p\le\iy$, invoking
the definition of the \HL\ maximal \g\ $\M(|\n f|)(x)$ (recalled below) we
get the following corollary.

\bco \label{c3}
For all $x\in\R^n$ the \iq y
\beq \label{e17}
MQf(x)\le \M(|\n f|)(x).
\e
 holds; in particular, $W^{1,p}(\R^n)\subset MQ^{1,p}(\R^n)$ for $1<p\le\iy$.
\eco

The inverse \iq y and inverse inclusion, though essentially true, are
somewhat more delicate \cite{BHS}, \cite{BKZ}, \cite{BBM2}, \cite{DL}, \cite{r5}, \cite{Ha1}.

Since the concept of the \HL\ maximal \f\ \cite{He} is very helpful
to describe
important consequences of \er{16} \cite{Gra}, \cite{HL}, we recall it briefly.

The centered \HL\ maximal \f\ of $f\in L^1\loc(\R^n)$ is defined as the supremum
of the \a es
\[
\M(f)(x)=\sup_{r>0} \fint{B(x,r)}{}|f|\,dy.
\]
The uncentered \HL\ maximal \f\ is defined as
\[
M(f)(x)=\sup_r \fint{x\in B(y,r)}{}|f|\,dz,
\]
i.e.\ when the \a ing is performed over all open balls that contain the
point~$x$.

Obviously $\M f(x)\le Mf(x)\le 2^n\M f(x)$ and both $\M f$ and $Mf$ are lower
semicontinuous, though, roughly speaking, $Mf$ is more regular (``smoother'')
than $\M f$, e.g.\ $Mf$ is even continuous.

For \f s on the real line, both
right $M_Rf$ and left $M_Lf$ are useful maximal \f s and play an important
role in many subtle questions of classical analysis on the real line $\R$:
\smallskip
\[
M_Rf(x)=\sup_{r>0}\hbox{\ooalign{\hfil$\relbar$\hfil\crcr$\displaystyle
\intop_{x}^{x+r}$}}|f(y)|\,dy, \quad
M_Lf(x)=\sup_{r>0}\hbox{\ooalign{\hfil$\relbar$\hfil\crcr$\displaystyle
\intop_{x-r}^x$}} |f(y)|\,dy.
\]
For \f s on $\R^n$ considered as the Cartesian products $\R^n=\R\times \ldots
\times\R$, $n$~times, also iterated maximal \f s \cite{JP}, \cite{Gra} and
\HL\ maximal \f s \wrt cubes, rectangular
maximal \f s, dyadic maximal \f s $M_df(x)$ are very useful in a variety of
contexts.

\bpr \label{p3}
For $1<p\le\iy$ the \pw\ \iq y \er9 characterises the \Ss\ $W^{1,p}(\R^n)$
and $W^{1,p}\loc(\R^n)$ $($and $MQ^{1,p}(\R^n)$ as well\/$)$.
\epr

\proof
Proposition \ref{p3} and its various proofs were well understood
(and published) at the time when \cite{BH} and \cite{Ha1} first appeared.
However, the proof of the \iq y \er9 for the class $MQ^{1,p}(\R^n)$
(Proposition~\ref{p1} above) is, apparently, published here
first\footnote{The present proof of Proposition~\ref{p1}, together with
Lemma~\ref{l1} above, also explicitly shows that the complete proof of
Lemma~1 in~\cite{Boj-gr} does not need any references to ``known facts of
harmonic analysis'' \cite{St1}, \cite{St2}. It is a direct and elementary consequence of the
accepted definitions of the spaces $MQ^{1,p}$ and the \Ss\ $W^{1,p}(\R^n)$.}.
We skip here the details of the well known proof that the \iq y
\er9 for a \f\ $f\in L^{p}\loc(\R^n)$ implies that $f\in W^{1,p}\loc(\R^n)$.
\endproof

Before recalling the rather not trivial implications of the \pw\ Sobolev \iq
y \er9,  the local \P\ type \pw\ inequalities \er{9.3}, \er{9.4}, \er{9.5} and the \P\
integral \iq ies \er{14}, \er{15} we show that the concept of the mean-\q\ spaces
$MQ^{1,p}\loc(\R^n)$ is meaningful in the general \mms\ context and, via
Proposition~\ref{p1} above, under some mild geometric conditions (doubling
measure and overlapping condition \er8, \er{9.2}) it leads directly to the metric
Sobolev--\P\ type vector spaces.

Thus, given a triple $(X,d,\mu)$: a complete \mms\ $X$ with distance
$d=d(x,y)$, $x,y\in X$ and a Borel regular measure~$\mu$, exactly as in \er1
and \er2 and the Definition~\ref{d1} above, we define the (restricted) maximal
mean \q\ ($d(z,x)=|z-x|$ in \er2)
\beq \label{e24}
MQ_R f(x)=\sup_{r<R} \fint{B(x,r)}{} \frac{|f(z)-f(x)|}{d(z,x)}\,d\mu
\e
and the maximal mean \df\ \q\ at $x\in X$
\beq\label{e25}
MQf(x)=\sup_{R>0} MQ_R f(x).
\e
The requirement $MQf\in L^p(X,\mu)$, $p\ge1$, where $L^p(X,\mu)$ denotes the
Lebesgue space of the measure space $(X,\mu)$, defines then the vector space
$MQ^{1,p}(X,d,\mu)$. If the geometric measure conditions, analogous to \er8
or \er{9.2} hold, the proof above of the Proposition~\ref{p1} with appropriate changes is valid.
We obtain

\bco \label{c4}
The linear vector space $MQ^{1,p}(X,d,\mu)$ has the structure of a metric
\Ss\ with intrinsic metric \g\ $MQf(x)\in L^p(X,\mu)$.
\eco

If the triple $(X,d,\mu)$ satisfies some kind of doubling condition or
Ahlfors $Q$-regular measure condition, the theory of Haj\l asz--Koskela and their followers applies
\cite{Ha1}, \cite{Ha2}, \cite{HaK1}, \cite{HaK2}, \cite{He}, \cite{HeK1},
\cite{HeK2}, \cite{Kos}.

In particular, the analogues of our Propositions \ref{p1}, \ref{p2}, \ref{p3}
above and the local \P-type \iq ies described above in Corollaries \ref{c1},
\ref{c2}, \ref{c3}, \ref{c4} hold.
These imply far reaching, mostly not trivial, consequences with proofs which
can be traced back to many authors.

We give below a short list of some of the results that can be deduced.
\ben
\item Refinements of the Lebesgue differentiation theorems for $L^1(X,\mu)$ \f s
for Sobolev $W^{1,p}(X,d,\mu)$ and $MQ^{1,p}(X,d,\mu)$ spaces. Approximate
differentiability in the Euclidean case.
\item H\"older continuity of \f s in $MQ^{1,p}(X,d,\mu)$ spaces for $p>s$,
Sobolev imbedding $MQ^{1,p}(X,d,\mu)\hookrightarrow C^\alpha\loc$, $p>s$,
$\alpha=1-\frac sp$, $\alpha>0$ ($s$---Hausdorff
dimension of $(X,d,\mu)$; in the Euclidean case $X=\R^n$, $s=n$).
\item Differentiability in the Euclidean case for $p>n$, A. P. Calderon theorem \cite{Cal},
J. Cheeger type approach in measure metric case for the triple $(X,d,\mu)$
\cite{BRZ}, \cite{Boj-Sb}, \cite{BIw}, \cite{Ch}, \cite{He}, \cite{Ke},
\cite{St1}.
\item Luzin type approximation by Lipschitz or $C^1$ \f s \cite{HdN}.
\item The Haj\l asz--Koskela Theorem \cite{HaK2} on the general local Sobolev imbedding
theorem of $MQ^{1,p}(X,d,\mu)$ spaces for $p< s$ into $L^{p^*}$ spaces,
$p^*=\frac{ps}{s-p}$.
\een

\section{Higher order Sobolev spaces}

For Sobolev space $W^{m,p}(\R^n)$, $m\ge 2$, $1\le p\le\infty$ higher degree
Taylor polynomials $T^kf(y;x)$, $0\le k\le m-1$, centered at $x$ and the
Taylor--Whitney remainders $R^kf(y;x)$ come into play.

Thus, by definition, for $f\in W^{m,p}(\R^n)$ we have
\begin{equation}\label{eq2.1}
T^kf(y,x)=f_0(x)+f_1(x)(y-x)+\ldots +f_k(x)\frac{(y-x)^k}{k!}
\end{equation}
and
\begin{equation}\label{eq2.2}
R^kf(y,x)=f(y)-T^kf(y,x)
\end{equation}
with
\[
f_0(x)\equiv f(x)
\]
and
\begin{equation}\label{eq2.3}
f_k(x)=D^k_xf(x),\ k=1,\dots,m-1
\end{equation}
understood as the generalized Sobolev, weak or distributional, derivatives of
the function $f$.

In the formulas (\ref{eq2.1})--(\ref{eq2.3}) we use standard polylinear
algebra and multidimensional analysis notations for the
$n$-dimensional vector space.

It is convenient to call the variables $x$ in the formulas
(\ref{eq2.1})--(\ref{eq2.3}) as \emph{field variables} and the variables $y$ as
\emph{space variables}: $y\in \R^n$. Although in our basic model case of Sobolev spaces:
$W^{m,p}(\R^n)$ both $x$ and $y$ vary in the whole space $\R^n$, $(x,y)\in
\R^n_x\times \R^n_y$, in most interesting cases, considered in the general
theory and in applications, the field variables are restricted to a subset
$\Sigma\subset \R^n$, a subdomain of $\R^n$, open or closed or just on
arbitrary closed subset of $\R^n_x$. Both variables $x$ and $y$ in the
Taylor--Whitney remainders $R^k f(y,x)$, and $R^k_l f(y,x)$ introduced later,
are considered to be field variables and take their values in the subset
$\Sigma$. For general Sobolev spaces the linear operator coefficients
$f_0,\dots,f_k$, $k\le m-1$ have their values in the Lebesgue spaces
$L^p(\R^n)$, $(L^p(\Sigma))$.

The operator valued coefficients $f_k$ in (\ref{eq2.1}) and (\ref{eq2.2})
need not satisfy the restricting conditions (\ref{eq2.3}). In this general
case we call the collection $F=\{f_0(x),\dots,f_k(x)\}$, $x\in\Sigma$ a
\emph{Whitney $k$-jet $F$}, defined on $\Sigma$ and the expressions (\ref{eq2.1})
and (\ref{eq2.2}) are termed the \emph{Taylor--Whitney fields} and \emph{Taylor--Whitney
remainders of the $k$-jet $F$}. For convenience the set of $k$-jets on
$\Sigma$ is denoted as: $J^k(\Sigma)$.

Taylor--Whitney fields $T^kF(y,x)$ can be differentiated in $y$, without any
restrictions and can be nicely described in terms of the formal jet
derivatives $D_l:J^k\to J^{k-|l|}$, $D_lF=\{f_l,f_{1+l},\dots,f_k\}$
\begin{equation}\label{eq2.4}
D^l_yT^kF(y,x)=T^{k-|l|}(D_lF)(y,x)
\end{equation}
and their Taylor--Whitney remainders
\begin{align}
&R^{k-|l|}(D_lF)(y,x)=R^k_lF(y,x),\  l=0,1,\dots,k-1,\label{eq2.5}\\
&R^kF(y,x)=R^k_0F(y,x) \nonumber
\end{align}
connected by the formulas (\ref{eq2.2}) for the $l$-th component of the
$k$-jet $F$ or the first (i.e.\ zero) component of $D_lF$
\begin{equation}\label{eq2.6}
R^{k-|l|}(D_lF)(y,x)+T^{k-|l|}(D_lF)(y,x)\equiv f_l(y).
\end{equation}
The recalled formulas (\ref{eq2.1})--(\ref{eq2.6}) and their admissible
derivatives with respect to the space variables, generate what has been called
the Taylor-algebra: $\T_m(\R^n)$. It has been introduced mainly by G.~Glaeser
\cite{Gla1} and B.~Malgrange \cite{Mal} as an algebraic tool to describe the
H.~Whitney \cite{Wh2}--\cite{Wh4} theory of continuously differentiable
functions $C^k(\Sigma)$, on arbitrary closed subsets of $\R^n$.

As shown in a series of papers, starting \cite{BH}, \cite{Boj5}, \cite{Boj13}
etc., up to the present one, these concepts, properly adjusted, are very useful
in the geometric understanding of the theory of Sobolev spaces as well.

Given a function $f\in W^{m,p}(\R^n)$ represented by its $(m-1)$ jet
\[
F=\{f(x),\dots,f^{(i)}(x),\dots,f^{(m-1)}(x)\}, \quad f^{(i)}=D^if
\]
we define the
$m$-th order restricted mean quotient
\begin{equation}\label{eq2.7}
M_RQ^m F(x)=\sup_{r\le R}\fint{B(x,r)}{}
\frac{\left|R^{m-1}F(z,x)\right|}{|z-x|^m}\,dz.
\end{equation}
We also define $m$-th maximal mean quotient
\begin{equation}\label{eq2.8}
MQ^m f(x)=\sup_{R>0}M_RQ^mf(x).
\end{equation}

In general the following fact holds:

\begin{prp}\label{prp2.1}
For a function $f\in W^{m-1,p}(\R^n)$, $1<p\le\infty$, the quantities
$MQ^m_Rf(x)$ completely control the Sobolev smoothness of the function $f$ in
$W^{m,p}(\R^n)$.
\end{prp}

As it stands for $m\ge 2$ Proposition~\ref{prp2.1} is rather vague and
needs some comments: The case $m=1$ is completely clarified by our
Proposition~\ref{p1} in Section~1:
we start with a function $f\in L^p(\R^n)$ and get
simple sufficient condition (obviously also necessary) for $f$ to be in the
Sobolev space $W^{1,p}(\R^n)$ with the full control of the norm\footnote{Here $L^{1,p}(\R^n)$
stands for the homogeneous space $W^{1,p}(\R^n)$ with the semi-norm
$\|f\|_{L^{1,p}}=\|\n f\|_{L^p}$.}
$||f||_{W^{1,p}(\R^n)}$: $MQ^1(f)\in L^p(\R^n)$,
$\|f\|_{L^{1,p}(\R^n)}\le C\|MQ^1f\|_{L^p}$ with an absolute constant $C$,
depending on~$n$ only.

For $m\ge 2$ the $m$-jet $F$ representing $f\in W^{m-1,p}(\R^n)$ allows us to
define the quantity $M_RQ^mf(x)$.

To extend Proposition~\ref{prp2.1} to arbitrary functions in the Lebesgue space
$L^p(\R^n)$ directly, some other process e.g.\ Newton--Lagrange interpolation,
described in \cite{Boj-gr}, see also \cite{Tri}, is needed. We shall come back
to this topic in more detail in \cite{16'}.

Now we briefly show how to handle the case $m=2$ with the elementary
geometric tools used for the case $m=1$ in Section~1 and some simple
considerations in the Taylor algebra $\T_2(\R^n)$, see \cite{Boj-f}.

\begin{lem}
For arbitrary $f\in W^{2,1}(\R^n)$, almost all $x,y\in \R^n$ and
$z\in\Sigma_r(x,y)$, $r=|x-y|$ the pointwise inequality holds
\begin{equation}\label{eq2.9}
\frac{\left|\R^{1}f(y,x)\right|}{|y-x|^2}\le
\frac{\left|\R^{1}f(z,y)\right|}{|z-y|^2}+
\frac{\left|\R^{1}f(z,x)\right|}{|z-x|^2}+
\frac{\left|f'(z)-f'(x)\right|}{|z-x|}.
\end{equation}
\end{lem}
This is an analogue of the estimate \er7 in Section~1.

\begin{proof}
To avoid the troublesome discussion of the set $N_f$ of  points in
$\R^n$, excluded by the words ``almost all'' above, we prove here the lemma for
$f\in C^1(\R^n)$. We have, $|N_f|=0$: the Lebesgue $n$-measure of
$N_f$ is zero. For $f\in C^1(\R^n)$, the proof reduces to simple algebraic
calculations in the Taylor algebra $\T_1(f)$ of order 1 and elementary
geometry of the sphere in $\R^n$ already exploited the discussion
of Sobolev spaces $W^{1,p}(\R^n)$.

For the triple of points $(x,y,z)\in \Sigma\times\Sigma\times\S$ we introduce
the expressions
\begin{align}
& Rf(y,x)=f(y)-f(x)\nonumber\\
&R^1f(y,x)=f(y)-f(x)-f'(x)(y-x)\label{eq2.10}
\end{align}
and analogous ones for the pairs $(x,z)$ and $(y,z)$.

{\makeatletter \c@eqn=\c@equation
\c@equation=0
\def\theequation{\thesection.\theeqn$_{\arabic{equation}}$}
\stepcounter{eqn}
In the Taylor algebra $\T_2(f)$ we have $(R\equiv R^0)$
\begin{gather}
R(y,x)+R(x,y)=0\label{e2.11.1}\\
R(y,x)-R(y,z)=-R(x,z)\label{e2.11.2}
\end{gather}
and in the Taylor algebra $\T_2(f)$ we have $(R\equiv R^0)$
\stepcounter{eqn}
\c@equation=0
\begin{align}
&R^1f(y,x)+R^1f(x,y)=\langle f'(y)- f'(x),y-x\rangle\label{e2.12.1}\\
P^1(x,y,z)&\equiv R^1f(y,x)-R^1f(y,z)\nonumber\\
&=-R^1f(z,x)+\langle f'(z)-f'(x),y-z \rangle \label{e2.12.2}
\end{align}
\c@equation=\c@eqn
\makeatother
}Now \er{2.12.2} for $z\in\Sigma_r(x,y)$, $r=|x-y|$ gives
\setcounter{equation}{12}
\begin{equation}\label{eq2.13}
\frac{\left|\R^{1}f(y,x)\right|}{|y-x|^2}\le
\frac{\left|\R^{1}f(z,y)\right|}{|z-y|^2}+
\frac{\left|\R^{1}f(z,x)\right|}{|z-x|^2}+
\frac{\left|f'(z)-f'(x)\right|}{|z-x|}.
\end{equation}
which completes the proof.
\end{proof}

It is also important to note the following consequence of \er{2.12.2}:
\begin{equation}\label{eq2.14}
D_yP^1(x,y,z)\equiv f'(z)-f'(x).
\end{equation}
Notice now that the integral means
\[
\fint B{}\frac{\left|P^1(x,y,z)\right|}{|z-y|^2}\,dz\quad \mbox{and}\quad
\fint B{}\frac{\left|P^1(x,y,z)\right|}{|x-z|^2}\,dy
\]
over balls $B$ in the space $\R^n_y$ are controlled by the values of
$MQ^2f(y)$ and $MQ^2f(x)$. In consequence, by the Markov inequalities, see
\cite{Boj13} and \cite{Boj-f}, the values of the quotients
$\frac{|f'(z)-f'(x)|}{|z-x|}$ are estimated from above by the maximal
functions of the gradient $\n^2f$ pointwise. Thus we see also that in this
somehow more subtle approach to pointwise inequalities for Sobolev functions
$f\in W^{m,p}(\R^n)$, $p>1$, $m\ge 2$, the iterated maximal functions of the
gradients appear naturally. The fundamental fact is that they all are, in the
sense of $L^p$ norms for $p>1$, estimated from above by $L^{p}$ norms of the Hardy--Littlewood
maximal functions $M(|\nabla^mf|)(x)$ or their iterates. In the existing
literature this seems not to have been exposed clearly enough.

Now, by averaging of (\ref{eq2.13}) over $\Sigma_r$ we obtain the pointwise
control of the left-hand side of (\ref{eq2.13}) by the mean maximal quotients
$MQ^2f$ and their Hardy-Littlewood maximal functions $M(MQ^2f)$ evaluated at
the points $x$ and $y$, as required in Proposition~\ref{prp2.1}.

We remark also that the averages
\[
\fint B{}\frac{|f'(z)-f'(x)|}{|z-x|}\,dz
\]
in (\ref{eq2.9}) can be also handled for functions $f$ in $W^{2,1}(\R^n)$, or
$f'\in W^{1,1}(\R^n)$ by an inductive procedure with respect to $m:f\in
W^{2,p}(\R^n)\to f'\in W^{1,p}(\R^n)$ and $f'\in MQ^{1,p}(\R^n)$, $p\ge 1$, that is
the resulting inclusion
\[
\fint B{}\frac{|f'(z)-f'(x)|}{|z-x|}\,dz\in L^p(\R^n)
\]
holds.

In this way we get yet another approach to the proof of the
Proposition~\ref{prp2.1}.


\section{Final comments}

Paragraph 2 above is the first introductory step in a novel approach to the
Sobolev pointwise inequalities characterizing the spaces $W^{m,p}(\R^{n})$
for $m\ge 2$. It will be continued in  \cite{16'}, along the road sketched in
the introduction. In particular we hope to analyze in \cite{16'} the
simplest, most natural and transparent proofs of the  inequality \eqref{0.25}
expressed in terms of the divided differences and finite differences, see
\cite{r0}, \cite{Boj-gr}, \cite{St2}.

Some of the results and ideas developed in \cite{BIK}, \cite{BKZ}
will be included in \cite{16'} as well.

We hope also to extend our approach to fractional Sobolev spaces, Besov type spaces and related trace and extension problems.

Stabilization to a polynomial of the functions in the homogeneous Sobolev spaces
$\mathop{W}\limits^{\circ}{}^{m,p}(\R^{n})$ will be also considered. See \cite{Ba}, \cite{BIN},
\cite{Wh5}.

Clearly, the remarks at the end of paragraph 1 and the list of problems quoted there,
should be considered as a forecast for a more detailed and thorough discussion in \cite{16'} based on
the results, methods and proofs scattered in  some of the references in our list at the end of this paper.

The proofs mentioned require a rather subtle refinement of the concept of Hardy--Littlewood
maximal function \cite{BIw}, \cite{Gra}, \cite{HaK2}, \cite{He}, \cite{HeK1}, \cite{Iw},
\cite{Ge}, \cite{HeKST}, \cite{Kos}, \cite{La}, upper gradient, higher order upper gradients,
weak maximal functions, weak and strong inverse H\"older
inequalities, etc. Here we refer to  \cite{BIw} where some of these refinements were elaborated
mainly in connection with the requirements of the analytical problems of quasiconformal mappings
theory, though some of them may be traced back even to some earlier works.

In the early eighties of the  last century the paper \cite{BIw} was rather widely studied in Finland
and, later, quite many of the concepts and methods used in \cite{BIw}, were
somehow ``pushed out'' to the
``folklore''!

They occupy a prominent role, though, in many papers of the Finnish School from
the period 1990--2005 quoted in our reference list.
As an example see the proof of A. Calder\'on theorem \cite{Cal} on total differentiability of
Sobolev functions in $W^{1,p}(\R^{n})$ for $p>n$, given in \cite{He}, originally in \cite{BIw}, see
also \cite{Boj-Sb}, \cite{GL}, \cite{Men},~\cite{St1}.

Sobolev space theory discussed along the novel approach propagated, among
others, in this paper and in \cite{BKZ}, suggests, and stimulates, asking
many new, apparently interesting and non-trivial questions. Here is an
example: in an ``admissible'' family $F$ of interpolation schemes $Z$,
invariant under affine transformations of the Euclidean space $\R^{n}$, for a
function $f$ in the Sobolev space $W^{m,p}(\R^{n})$, describe an
``optimal''(?) one, ``approximately optimal''(?!), for its ``numerical
efficiency''(?) or ``cost of its computational
implementation''(?)\footnote{question marks mean that the suggested concepts
{\it are to be defined as an essential part of the research process},
requiring maybe a rather deep, new insight.}

As a working conjecture for this problem we state the following
inequality for the Sobolev smoothness density $g_{f,Z}(x)$
\beq\label{3.1}
g_{f,Z}(x) \le C_{n,m}g_{f,TW}(x),
\e
where $g_{f,TW}(x)$ is the notation of the Sobolev smoothness of the function
$f$ for the Taylor--Whitney interpolation, as in \eqref{eq2.2} or in \cite{BH}, and $g_{f,Z}$
is calculated in analogy with \eqref{0.3}, see also \cite{BKZ};
 $C_{n,m}$ is an absolute constant to be found.

Even for the case of Sobolev functions
on the real line the inequalities of type \eqref{3.1} and the constant in \eqref{3.1} are
not known \cite{DL}, \cite{r8}, \cite{BKZ}.

In the ``enormous'' literature on Sobolev spaces and their applications there can
hardly be found any references concerning this and related questions --- even
on a segment of the real line! --- \cite{4'}, \cite{BIw},
\cite{BBM2}, \cite{r4}.

Referring again to the ``novelty'' of the approach to the general theory of
Sobolev spaces, emerging from this paper and others recent related, let me express the hope that
this will be a useful contribution in the direction of producing a concise, readable exposition
of the Sobolev theory, with minimal necessary technicalities\footnote{Man
soll die Dinge so einfach machen wie m\"oglich --- aber nicht einfacher (A. Einstein);
translation for sciences (H.~Triebel): One should present assertions
as simple as possible --- but not simpler.}
and not missing
any typical properties and concepts currently used in geometrical, analytical and
interdisciplinary applications.
Any critical remarks, suggestions, corrections and new references will be very welcome.

Despite the rapidly growing list of recent monographic publications\footnote{A. and Yu. Brudnyi's I, II, 2011, V. G. Mazya 2011, W. Yuan, W. Sickel, D. Yang 2010, S. Kislyakov, N. Kruglyak 2013,
L. Pick, A. Kufner, O. John, S. Fucik I, II 2013.} on topics of Sobolev spaces theory some farther
effort seems to be worthwhile.




\end{document}